\documentclass{amsart}
\pdfoutput=1
\usepackage{graphicx}

\newtheorem{thm}{Theorem}[section]    % Standard theorem environment
\newtheorem{lem}[thm]{Lemma}          % Lemma environment with numbering 
\newtheorem{cor}[thm]{Corollary}
\newtheorem{prop}[thm]{Proposition}
\newtheorem*{thmmainexample}{Theorem~\ref{mainexample}}
\newtheorem*{thmsbtheorem}{Theorem~\ref{sbtheorem}}
\newtheorem*{corbound}{Corollary~\ref{bound}}
\newtheorem*{propwhcalc}{Proposition~\ref{whiteheadcalculation}}
\theoremstyle{definition}
\newtheorem{defn}[thm]{Definition}    % Definition environment with 
\newtheorem{exam}[thm]{Example}
\newtheorem*{rem}{Remark}             % Unnumbered environment for remarks.

\newcommand{\inv}{^{-1}}
\newcommand{\p}{\pi_1}
\newcommand{\del}{\partial}
\renewcommand{\a}{\alpha}
\renewcommand{\b}{\beta}
\newcommand{\lk}{{\mathrm{lk}}}

\newcommand{\SR}{{\mathcal{R}}}
\newcommand{\R}{{\mathbb{R}}}

\title{The First-order Genus of a Knot}

\author{Peter Horn}
\begin{document}

\begin{abstract}    % type your abstract below
	We introduce a geometric invariant of knots in $S^3$, called the first-order genus, that is derived from certain 2-complexes called gropes, and we show it is computable for many examples.  While computing this invariant, we draw some interesting conclusions about the structure of a general Seifert surface for some knots.
\end{abstract}

\maketitle

%%%%%%%%%%%%%%%%%%%%   Start of main body of article

\section{Introduction}

The main objective of this paper is to define and investigate a geometric knot invariant that is stronger than the genus of a knot and that is in fact a ``higher-order genus" of a knot.  The invariant is called the \emph{first-order genus} and is denoted $g_1$.  Roughly speaking, one obtains the first-order genus of $K$ by adding the genera of certain curves on a minimal genus Seifert surface for $K$.  Further details will be provided in Section~\ref{definitions}.  We prove that the first-order genus is bounded below by the three-dimensional genus (cf. Proposition~\ref{twice}) and hence is more suitable than the slice genus for distinguishing infinitely many knots that have the same genus.  We also prove that our invariant is independent of many classical three- and four-dimensional geometric and algebraic knot invariants.

\begin{thmmainexample}
	There exists an infinite family of (distinct) knots $\mathcal{H}$ with the property that for all $K, J\in\mathcal{H}$ with $K\neq J$, the following hold:
	\begin{enumerate}
		\item $g(K)=1$,
		\item $\Delta_K(t)=1$,
		\item $K$ is smoothly slice, and
		\item $g_1(K)\neq g_1(J)$.
	\end{enumerate}
\end{thmmainexample}

To compute the first-order genus of $K$, one must inspect all minimal genus Seifert surfaces of $K$ (details are again postponed until Section~\ref{definitions}).  Thus, any information about an arbitrary minimal genus Seifert surface helps in the computation of the first-order genus.  For knots $J$ and $L$ and integers $m$ and $n$, let $K(J,L,m,n)$ denote the knot pictured below:

\begin{center}
	\begin{picture}(140, 100)(0,0)
			\includegraphics[scale=.8]{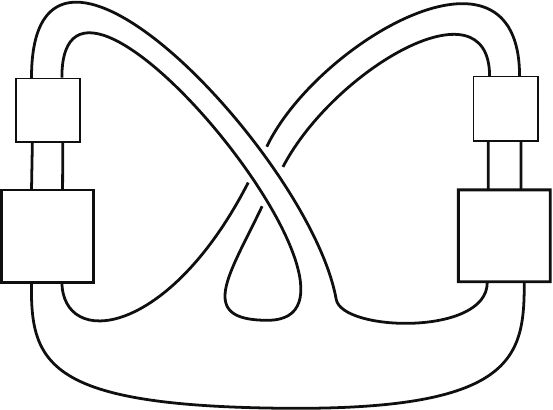}
			\put(-121, 67){\small $m$}
			\put(-13, 67){\small $n$}
			\put(-121, 37){$J$}
			\put(-14, 37){$L$}
		\end{picture}
\end{center}

Using cut-and-paste techniques common to three-dimensional topology, we are able to prove the following, which tells us a great deal of information about any minimal genus Seifert surface for $K(J,L,m,n)$.

\begin{thmsbtheorem}
	If $J$ and $L$ are neither trivial nor cable knots, then any minimal genus Seifert surface for $K(J,L,m,n)$ has a symplectic basis $(\alpha,\beta)$ where $\alpha$ and $J$ have the same knot type, and $\beta$ and $L$ have the same knot type.
\end{thmsbtheorem}

\begin{corbound}
	Under the hypotheses of Theorem~\ref{sbtheorem}, $g_1(K(J,L,m,n))\geq g(J)+g(L)$.
\end{corbound}

Now consider the following diagram of the unknot (on the left):

\begin{figure}[ht]
	\begin{center}
	\begin{picture}(195, 55)(0,0)
		\put(94, 25){$\to$}
		\put(118, 28){\small $K$}
		\includegraphics{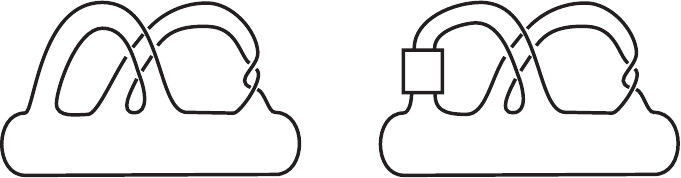}
	\end{picture}
	\end{center}
\end{figure}

One can construct a new knot by taking the left-hand band (which consists of two strands) and tying it into a knot $K$.  This is the well-known Whitehead double of $K$, denoted $Wh_0(K)$ and depicted on the right, which has genus one regardless of the knot $K$.  If two knots have the same genus, one might try to distinguish the knots by showing they have different slice genera.  This approach might be fruitful when trying to distinguish a finite number of knots of the same genus but will fail when investigating any infinite family.  While the knot genus and topological slice genus do not distinguish untwisted Whitehead doubles, the following proposition indicates that the first-order genus will sometimes distinguish them.

\begin{propwhcalc}
      Let $J$ be a nontrivial and noncable knot, and let $Wh_0(J)$ denote the (positively- or negatively-clasped) untwisted Whitehead double of $J$. Then
      $g_1(Wh_0(J))=1+g(J)$.
\end{propwhcalc}

We use our invariant to distingish Whitehead doubles not studied by Brittenham and Jensen in~\cite{BJ06}.

In the early 1950s, Schubert found a lower bound for the genus of a satellite knot~\cite{hS53}.  More specifically, if $K$ is a winding number $n$ satellite of $\widehat K$, where $\widetilde K$ is the pattern knot, we have $$g(K)\geq n\,g\left(\widehat K\right) +g\left(\widetilde K\right)$$
If $K$ is a winding number zero satellite of $\widehat K$, Schubert's inequality sometimes fails to give an accurate lower bound for $g(K)$.  For example, think of the Whitehead double of a knot to be a winding number zero satellite whose pattern is the unknot.  Our Proposition~\ref{whiteheadcalculation} and Corollary~\ref{bound} are evidence that our invariant is more suited for distinguishing winding number zero satellites.

Cochran~\cite{tC04} defined the higher-order Alexander modules and higher-order linking forms, and Cochran-Orr-Teichner~\cite{COT03} introduced higher-order $L^2$-signatures to the study of knot concordance.  As demonstrated in ~\cite{tC04} and ~\cite{COT03}, these higher-order algebraic invariants can distinguish knots with identical classical algebraic invariants, provided the classical Alexander module is nontrivial.  In contrast, our invariant does distinguish some knots with trivial Alexander module (cf. Theorem~\ref{mainexample}).  To the author's knowledge, the first-order genus is the first higher-order \emph{geometric} invariant.

\section{Definitions}\label{definitions}

Let $\Sigma$ be a once-punctured, orientable surface of genus $g>0$.  A \textbf{symplectic basis of curves} for $\Sigma$ is a collection $\a_1,\b_1,\ldots,\a_g,\b_g$ of simple closed curves that satisfies the following:
\begin{enumerate}
	\item $\a_1,\b_1,\ldots,\a_g,\b_g$ forms a basis for $H_1(\Sigma)$,
	\item $\a_i\cap\b_j$ is a point if $i=j$ and is empty if $i\neq j$,
	\item $\a_i\cap\a_j$ is empty if $i\neq j$, and
	\item $\b_i\cap\b_j$ is empty if $i\neq j$.
\end{enumerate}

If $\Sigma$ has positive genus, then $\Sigma$ has infinitely many symplectic bases.  We now recall the definition of a grope.  A \textbf{grope} is a 2-complex that is formed by gluing once-punctured, orientable surfaces (henceforth ``surfaces") in stages.  A grope of height one is just a surface $\Sigma$.  Given a symplectic basis $\a_1,\b_1,\ldots,\a_g,\b_g$ for $\Sigma$, we construct a grope of height $h+1$ by attaching gropes of height $h$ to each $\a_i$ and $\b_i$ along the boundary circles (of the height $h$ gropes). Teichner~\cite{pT04} has a wonderful description of different types of gropes.  Gropes have appeared recently in filtrations of the knot concordance group (cf. ~\cite{COT03} and~\cite{CT07}).

If $K$ is a knot in $S^3$ with Seifert surface $\Sigma$ and $J$ is a simple closed curve on $\Sigma$, then $\lk(K,J)=0$; that is, $J$ is nullhomologous in $S^3-K$.  Thus, we can find an orientable surface embedded in $S^3-K$ that is bounded by $J$.  Define $g(J;K)$ to be the minimum of the genera of orientable surfaces embedded in $S^3-K$ that are bounded by $J$.  If $F: S^3\times[0,1]\to S^3$ is an isotopy, then $g(F(J,0);F(K,0)=g(F(J,1);F(K,1))$.

\begin{defn}\label{fogdef}
	Define the first-order genus of the unknot to be zero.  If $K$ is a nontrivial knot, define the \textbf{first-order genus} of $K$ to be $$g_1(K)=\min\left\{\min\left\{g(\a_1;K)+g(\b_1;K)+\cdots+g(\a_g;K)+g(\b_g;K)\right\}\right\}$$ where the innermost minimum is taken over all symplectic bases of a given minimal genus Seifert surface, and the outermost minimum is taken over all minimal genus Seifert surfaces for $K$.
\end{defn}

\begin{rem}
	It is important in the proofs of some results below that the (first stage) Seifert surface be of minimal genus.
\end{rem}

Since each $g(\a_i;K)$ is invariant under isotopies of $S^3$, the first-order genus of $K$ is an isotopy invariant of $K$.  We explain presently that $g_1(K)$ measures the geometric complexity of a certain mapped-in grope of height two that is bounded by $K$.  The first stage of this grope is a minimal genus Seifert surface for $K$.  To each $\a_i$ (respectively $\b_i$) on $\Sigma$, attach a surface with genus $g(\a_i;K)$ (respectively $g(\b_i;K)$); these second-stage surfaces are embedded in $S^3-K$ with boundary $\a_i$ (or $\b_i$) and may intersect (in their interiors) the first-stage surface $\Sigma$ and the other second-stage surfaces.  We call this type of grope a \textbf{weak grope} of height two.  The first-order genus of $K$ is the smallest sum of the genera of the second-stage surfaces of a weak grope of height two that is bounded by $K$.

As described above, any knot in $S^3$ bounds a weak grope of height two.  This geometric fact is slightly stronger than the algebraic fact that the longitude of $K$ lies in the second derived subgroup of $\p\left(S^3-K\right)$.  In fact, a knot $K$ bounds a mapped-in grope of height $n$ if and only if the longitude lies in the $n$-th derived subgroup of $\p\left(S^3-K\right)$ (cf. Teichner~\cite{pT04}).  One could define the $n$-th-order genus of a knot in a fashion similar to our Definition~\ref{fogdef}, but such a definition would be valid only for knots whose longitude lies in the $n$-th derived subgroup of $\p\left(S^3-K\right)$.  This is a farily big restriction, as~\cite[Proposition 12.5]{tC04} states that if the preferred longitude lies in the third derived subgroup, then the Alexander polynomial is trivial.

We now state some basic results.

\begin{prop}[Injectivity]\label{injectivity}
   $g_1(K)=0$ if and only if $K$ is the unknot.
\end{prop}
\begin{proof}
   Suppose $K$ is a knot with $g_1(K)=0$ and $g(K)=g>0$.  Pick a Seifert surface $\Sigma$ and symplectic basis
   $\a_1,\b_1,\ldots,\a_g,\b_g$ that realize $g_1(K)=0$.  Let $x$ be a point on $K$ and for $i=1,\ldots,g$, let $\gamma_i$ denote a path from $x$ to the point $\a_i\cap\b_i$.  Since $g_1(K)=0$,
   each $\a_i$ and $\b_i$ bounds a disc in $S^3-K$, hence
   $\gamma_i\a_i\overline{\gamma_i}=\gamma_i\b_i\overline{\gamma_i}=1\in\p\left(S^3-K,x\right)$ for each $i$.  Let $\ell$ denote the preferred longitude of $K$.  Since
   $\ell=\gamma_1[\a_1,\b_1]\overline{\gamma_1}\cdots\gamma_g[\a_g,\b_g]\overline{\gamma_g}$, $\ell=1\in\p\left(S^3-K\right)$.  By the
      unknotting theorem, $K$ is the unknot.
\end{proof}

\begin{prop}\label{twice}
 For any knot $K$, $g_1(K)\geq 2g(K)$.
\end{prop}
\begin{proof}
   If $K$ is the unknot, we have nothing to prove.  Otherwise,  let $\alpha$ be any basis curve for $H_1(\Sigma)$, where $\Sigma$ is a minimal genus Seifert surface for $K$.  If $\alpha$ bounds a disc embedded in $S^3-K$, then $\ker\left(\p(\Sigma)\to\p\left(S^3-K\right)\right)\neq 1$.  By a well-known lemma in 3-manifold topology \cite[Lemma 6.1]{jH76}, there is a disc $D\hookrightarrow S^3-K$ with $D\cap \Sigma=\partial D$ and $[\partial D]\neq 1\in\p(\Sigma)$.  Let $\beta_0=\partial D\subset \Sigma$, and let $\beta_1$ denote a small pushoff of $\beta_0$ in $S$, i.e. $\beta_0\coprod\beta_1$ bounds an annulus $A$ embedded in $\Sigma$.  Observe $\beta_0\cap \beta_1=\emptyset$ and $D\cap \Sigma=\beta_0$ imply $D\cap \beta_1=\emptyset$; thus ${\mathrm{lk}}(\beta_0,\beta_1)=0$, and the annulus $A$ extends to an embedding $D^2\times [0,1]\hookrightarrow S^3-K$ such that $D^2\times\{i\}=\beta_i$ for $i=0,1$, $A=\partial D^2\times I$, and $\left(D^2\times I\right)\cap \Sigma=A$.  Therefore $\left(\Sigma-A\right)\cup D_0\cup D_1$ is an embedded Seifert surface for $K$ with genus strictly smaller than that of $\Sigma$, contradicting $g(\Sigma)=g(K)$.
\end{proof}

\begin{cor}
	There exist knots with arbitrarily high first-order genera.
\end{cor}

Theorem~\ref{mainexample} complements the preceding corollary by proving the existence of genus one knots with arbitrarily high first-order genera.

\begin{prop}[Subadditivity]\label{subadditivity}
	For any knots $K$ and $J$, we have $g_1(K\# J)\leq g_1(K)+g_1(J)$.
\end{prop}
\begin{proof}
	We omit the technical details of this proof but give the main idea.  Any point on $K$ has a three-dimensional neighborhood that does not intersect the second-stage surfaces of a weak grope realizing $g_1(K)$.  We can pick any point on $K$ and any point on $J$ and perform the connected sum along these small neighborhoods.  The result is a weak grope of height two that bounds $K\# J$, and its first-stage surface is of minimal genus.  Adding up the genera of the second-stage surfaces yields the desired inequality.
\end{proof}

\section{Examples}\label{examples}

In this section, we supply many examples in hopes of convincing the reader that our invariant is computable in many cases.

\begin{exam}
	Consider the right-handed trefoil $T$.  Below is a picture of the (unique) minimal genus Seifert surface for $T$ with a symplectic basis of curves.
\begin{center}
	\begin{picture}(80, 85)(0,0)
		\includegraphics{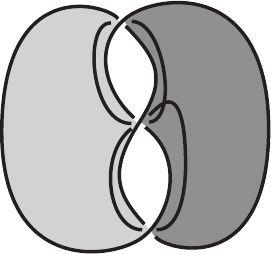}
	\end{picture}
\end{center}

One sees that each of the basis curves bounds a punctured torus that does not intersect $T$ by pushing each curve off of the surface and then tubing around part of $T$.  The sum of the genera of the basis curves is equal to $2$, and in light of Proposition~\ref{twice}, we conclude $g_1(T)=2$.
\end{exam}

A similar argument will prove that the first-order genus of the figure-eight knot is also equal to $2$.  One might ask if the first-order genus can distinguish knots of genus one.  We answer in the affirmative.

\begin{prop}\label{whiteheadcalculation}
      Let $J$ be a nontrivial and noncable knot, and let $Wh_0(J)$ denote the (positively- or negatively-clasped) untwisted Whitehead double of $J$. Then
      $g_1(Wh_0(J))=1+g(J)$.
\end{prop}
\begin{proof}
      Let $J$ be a knot which is neither the unknot nor a cable knot.  Let $K=Wh_0(J)$ be
      the untwisted Whitehead double of $J$.  Consider the genus one
      Seifert surface $\Sigma$ with symplectic basis $\a$ and $\b$
      depicted in Figure~\ref{whiteheaddiagram} ($\a$ is the labeled curve, and $\b$ is not labeled).

\begin{figure}[ht]
	\begin{center}
		\begin{picture}(350, 125)(0, 0)
			\includegraphics[scale=.6]{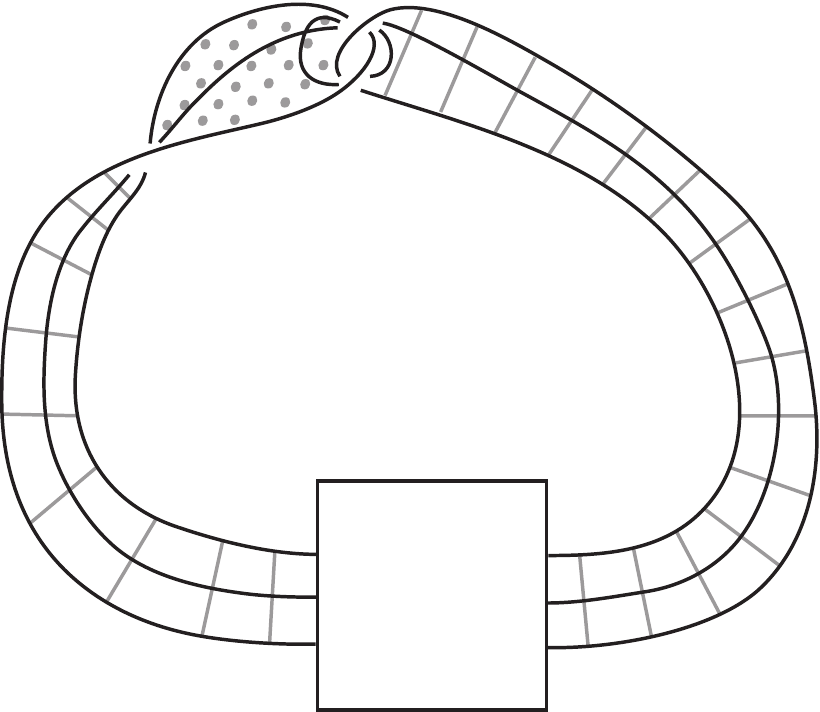} 
			\hspace{1cm}  \includegraphics[scale=.6]{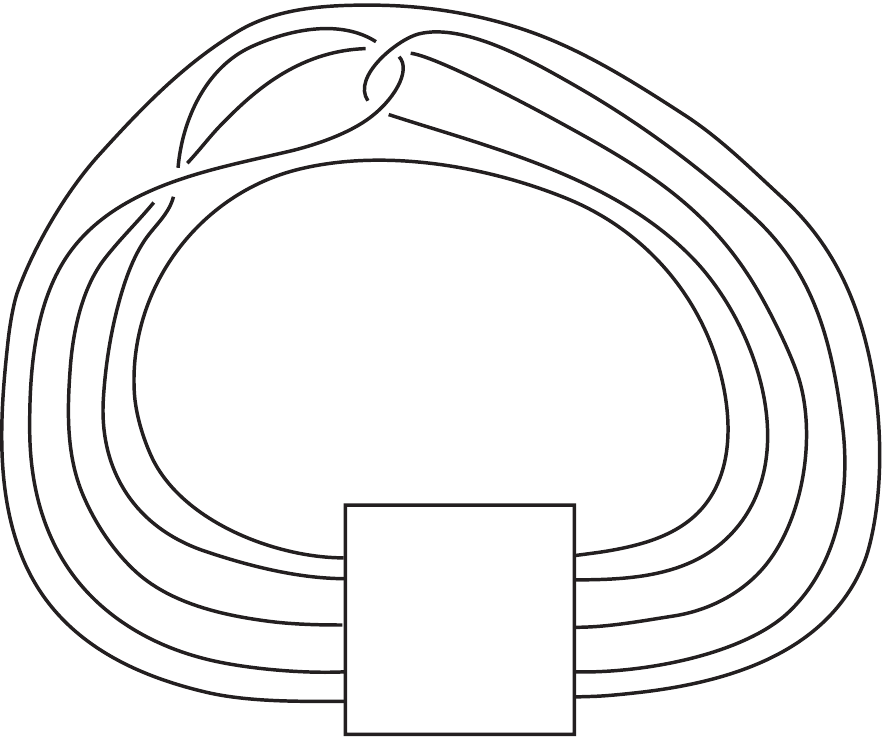}
			\put(-258, 16){$J$}
			\put(-304, 18){\small $\a$}
			\put(-214, 4){$K=\del \Sigma$ }
			\put(-76, 16){$J$}
			\put(-39, -2){$N(J)$}
			\put(-114, 23){\small$\a$ }
		\end{picture}
	\end{center}
	\caption{The seifert surface with the basis curves, and the knot inside the companion torus}
	\label{whiteheaddiagram}
\end{figure}
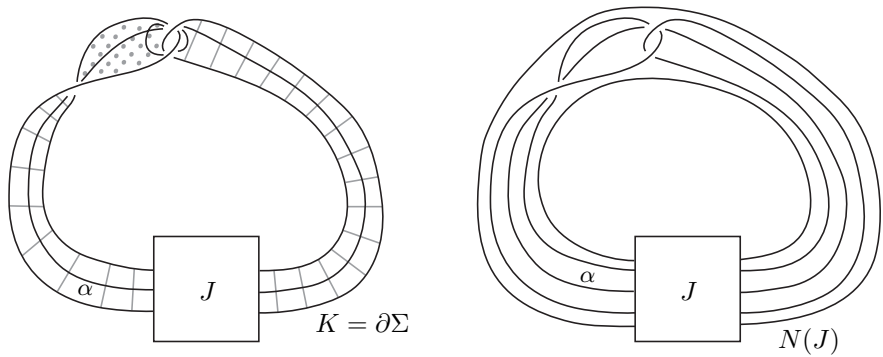

      One can see that $\b$ bounds an imbedded, punctured torus by
      tubing around $K$.  Notice that $\a$ has genus $g(J)$, so
      any imbedded surface
      in $S^3-K$ with boundary $\a$ must have genus at least $g(J)$.
      Pushing $\a$ off $\Sigma$ in the downward direction yields a
      parallel copy $\a^-$ of $J$ that does not link the knot $K$ (recall
      $K$ is untwisted); in fact, $K\cup \a^-$ is a boundary link
      with $K=\partial \Sigma$ and $\a^-$ the boundary of a (parallel copy of a) minimal genus
      Seifert surface for $J$.  Thus, $\a$ bounds an imbedded,
      oriented surface of genus $g(J)$, and $g(\a;K)+g(\b;K)=1+g(J)$.

      Now let $x$ and $y$ be any symplectic basis for $\Sigma$.  Without loss of generality,
      suppose $x=n\a+m\b$ where $n> 0$.  We see $x$ is a
      satellite of $J$ by Figure~\ref{whiteheaddiagram}.

      By a result of Schubert \cite{hS53}, $g(x)\geq n\, g(J)$, where $g(x)$ is the
      genus of $x$ in $S^3$.  Since $g(x;K)$ is at
      least $g(x)$ and $g(y;K)$ is at least one (by Proposition~\ref{twice}), we
      see $g(x;K)+g(y;K)\geq g(J)+1$.  Since $\{x,y\}$ is an arbitrary
      symplectic basis for $\Sigma$, $\min\{g(x;K)+g(y;K)\}\geq
      g(J)+1$, where the minimum is taken over all symplectic bases
      of $\Sigma$.

      A result of Whitten~\cite{wW73} guarantees that any minimal genus
      Seifert surface for $K=Wh_0(J)$ is isotopic to $\Sigma$, provided $J$ is a nontrivial and noncable knot.
      As this is the case at hand, $g_1(K)\geq g(J)+1$.  Since we realized a triple $(\Sigma,\a,\b)$ with
      $g(\a;K)+g(\b;K)=g(J)+1$, we conclude $g_1(K)=g_1(Wh_0(J))=g(J)+1$.
\end{proof}

\begin{prop}\label{pretzel}
	Let $K_n$ be the knot bounding the Seifert surface $V_n$ pictured below:
	\begin{center}
		\begin{picture}(240, 115)(0,0)
			\includegraphics[scale=.6]{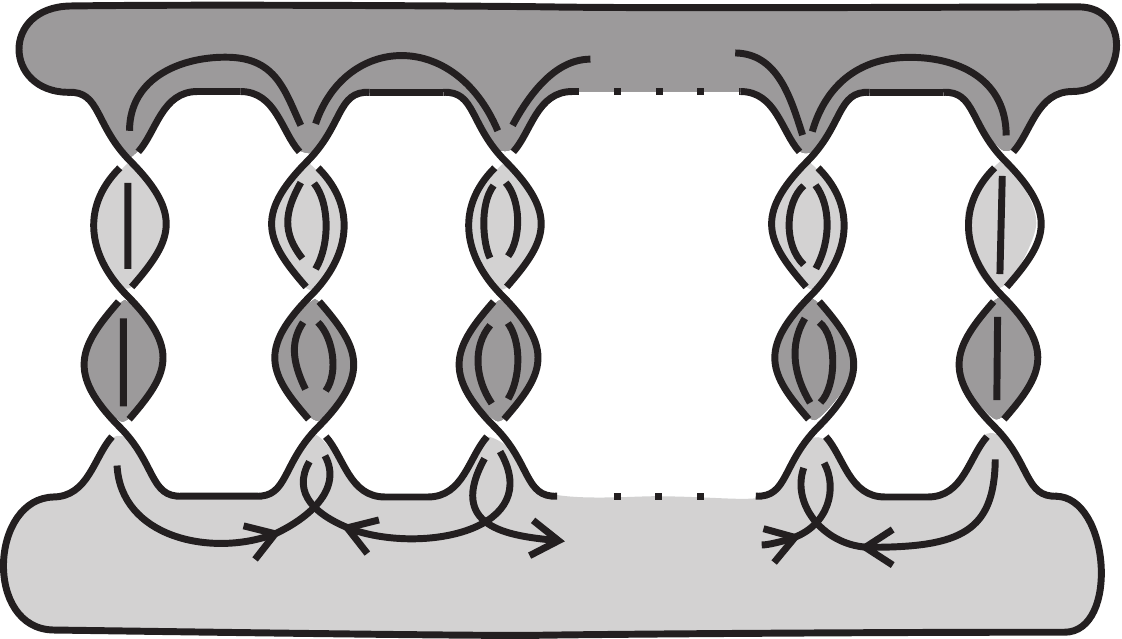}
			\put(5, 8){\small $K_n=\del V_n$}
			\put(-164, 11){\small $x_1$}
			\put(-129, 11){\small $y_1$}
			\put(-109, 11){\small $x_2$}
			\put(-58, 11){\small $x_n$}
			\put(-39, 11){\small $y_n$}
			\put(-189, 11){\small $+$}
			\put(-187, 100){\small $-$}
			\put(-94, 66){\small $2n+1$}
			\put(-94, 56){\small twisted}
			\put(-94, 46){\small strands}
		\end{picture}
	\end{center}

	Then for each $n\geq 1$, $g(K_n)=n$, $K_n$ is a not a cable knot, and $K_n$ is a ribbon knot.
\end{prop}
\begin{proof}
	Following the notation of Kawauchi~\cite{aK84}, $K_n=K(-3,3,-3,\ldots,3,-3)$, where the final `$-3$' is the $(2n+1)$-th entry.  The $K_n$ are pretzel knots and hence simple by~\cite{aK84}.  The genus of the surface $V_n$ is $n$, and $x_1,y_1,\ldots,x_n,y_n$ form a (non-symplectic) basis for $V_n$.  One computes the Seifert matrix of $V_n$ to be

$$\theta_n=\left(\begin{array}{ccccccccc}
      0&2&0&&&&&&\\
      1&0&-1&0&&&&&\\
      0&-2&0&2&0&&&&\\
      &0&1&0&-1&0&&&\\
      &&\ddots&\ddots&\ddots&\ddots&\ddots&&\\
      &&&0&-2&0&2&0&\\
      &&&&0&1&0&-1&0\\
      &&&&&0&-2&0&2\\
      &&&&&&0&1&0
\end{array}\right)$$

and one expands by the first column of $\theta_n-t\theta_n^T$ to prove inductively that
$\displaystyle \Delta_{K_n}(t)=\det\left(\theta_n-t\theta_n^T\right)=\left(-2t^2+5t-2\right)^n$.  We conclude that $V_n$ is a minimal genus Seifert surface for $K_n$ and $g(K_n)=n$.

We see a ribbon disc for $K_2$:
	\begin{center}
		\begin{picture}(200, 94)(0,0)
			\includegraphics[scale=.5]{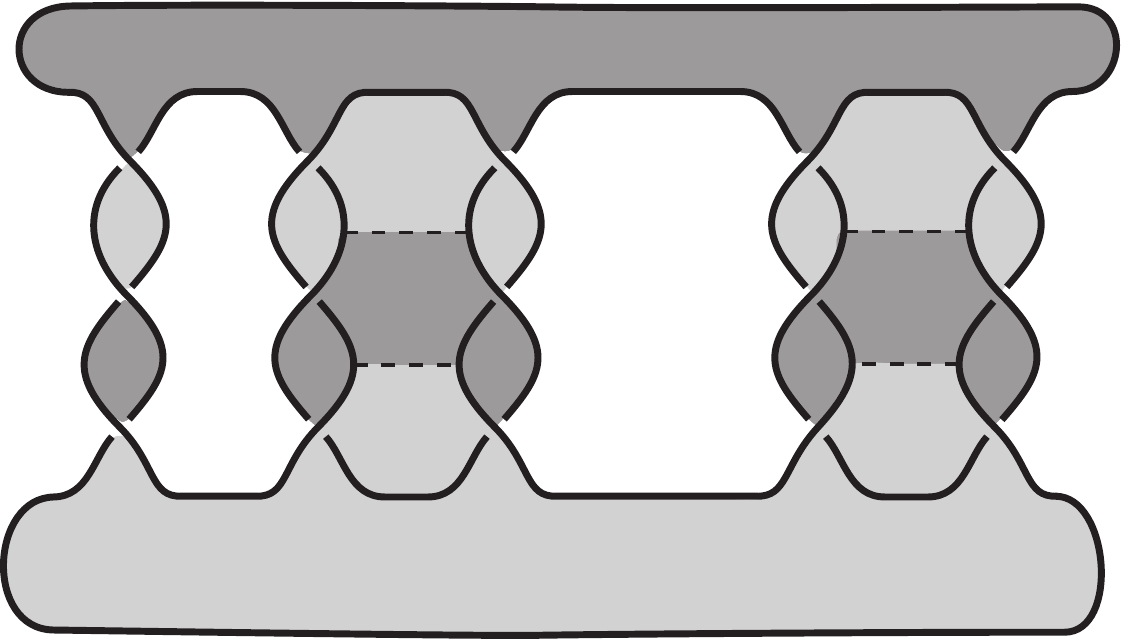}
		\end{picture}
	\end{center}

A ribbon disc for $K_n$ may be drawn in a similar fashion.
\end{proof}

We now have infinitely many knots with identical classical invariants but distinct first-order genera.

\begin{thm}\label{mainexample}
	There exists an infinite family of (distinct) knots $\mathcal{H}$ with the property that for all $K, J\in\mathcal{H}$ with $K\neq J$, the following hold:
	\begin{enumerate}
		\item $g(K)=1$,
		\item $\Delta_K(t)=1$,
		\item $K$ is smoothly slice, and
		\item $g_1(K)\neq g_1(J)$.
	\end{enumerate}
\end{thm}
\begin{proof}
	Let $\mathcal{H}=\left\{Wh_0(K_n):n\geq 1\right\}$ where $K_n$ are defined as in Proposition~\ref{pretzel}.  It is well-known that the untwisted Whitehead double of a smoothly slice knot is smoothly slice and that the Alexander polynomial of such a knot is trivial.  Proposition~\ref{whiteheadcalculation} implies $g_1(Wh_0(K_n))=1+n$.
\end{proof}

\begin{rem}
	In~\cite{BJ06} Brittenham and Jensen proved the canonical genus of a Whitehead double of a certain type of pretzel knot $K$ was equal to the crossing number of $K$.  In Kawauchi's notation, the pretzel knots of Brittenham and Jensen were of the form $K(k_1,\ldots,k_n)$, where $k_1,\ldots,k_n\geq 1$.  Since our pretzel knots are of the form $K(-3,3,-3,\ldots,3,-3)$, our first-order genus can distinguish knots that were not studied in~\cite{BJ06}.
\end{rem}

\section{First-order Genus and Seifert Surfaces}\label{fogandss}

We saw in Proposition~\ref {whiteheadcalculation} that if a knot has a unique minimal genus Seifert surface, the calculation of the knot's first-order genus is greatly simplified.  In this section, we calculate the first-order genus of a large family of genus one knots by proving that any minimal genus Seifert-surface for our knots must have a particular symplectic basis.  As in Proposition~\ref{whiteheadcalculation}, our present knots are satellites.

\begin{figure}[ht]
	\begin{center}
		\begin{picture}(140, 100)(0,0)
			\includegraphics[scale=.8]{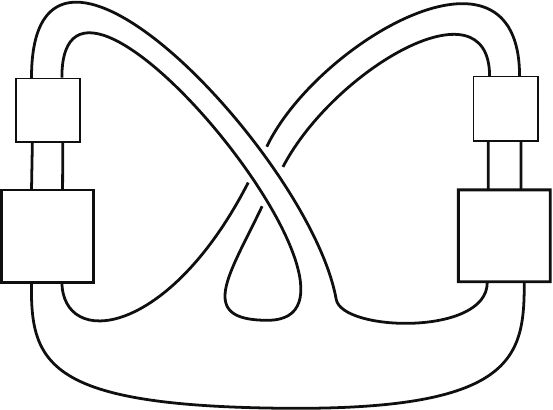}
			\put(-121, 67){\small $m$}
			\put(-14, 67){\small $n$}
			\put(-121, 38){$J$}
			\put(-14, 38){$L$}
		\end{picture}
		\caption{The knot $K(J,L,m,n)$}
		\label{knot}
	\end{center}
\end{figure}

\begin{defn}
	For knots $J$ and $L$ and integers $m$ and $n$, define $K(J,L,m,n)$ to be the knot depicted in Figure~\ref{knot}; the left-hand strands are tied into the knot $J$ and twist $m$ full times around each other.  When $m=n=0$, we denote this knot by $K(J,L)$.  In addition, let $\SR$ denote the set of nontrivial knots in $S^3$ which are neither torus nor cable knots.
\end{defn}

The purpose of this section is to obtain a lower bound for the first-order genus of $K(J,L,m,n)$ when $J,L\in\SR$.  In finding this lower bound, we discover a remarkable fact about any minimal genus Seifert surface for $K(J,L,m,n)$ (cf. Theorem~\ref{sbtheorem}).

\subsection{Certain $K(J,L,m,n)$ are Satellite Knots}

\begin{figure}[ht]
  \begin{center}
  	\begin{picture}(260, 90)(0,0)
     	\includegraphics[viewport=0 0 255 85, scale=.8]{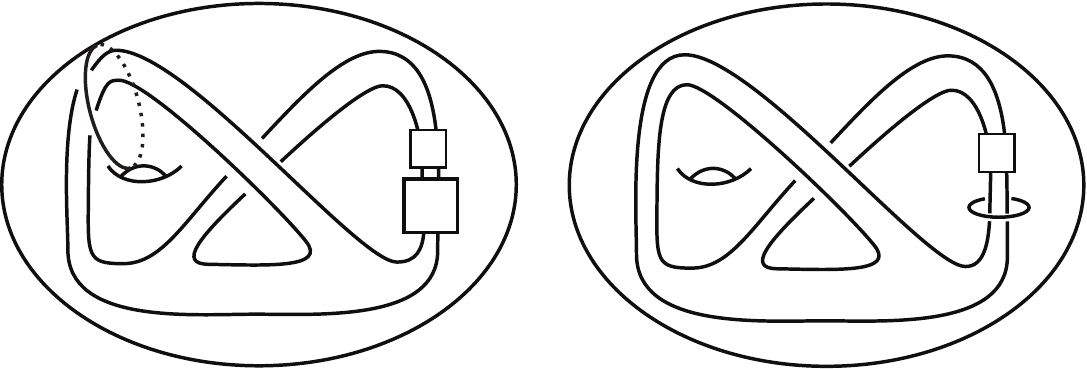}
		\put(-199, 7){\small $V$}
		\put(-66, 7){\small $V$}
		\put(-108, 35){\small $L$}
		\put(-108, 49){\small $n$}
		\put(-164, 69){\small $P_{L,n}$}
		\put(-33, 69){\small $P_{n}$}
		\put(-192, 77){\small $\mu$}
		\put(23.5, 48){\small $n$}
		\put(34, 36){\small $\eta$}
	\end{picture}
  \end{center}
  \caption{The pattern knots $P_{L,n}$ and $P_n$, as well as the solid torus $\eta$}
  \label{pattern}
\end{figure}
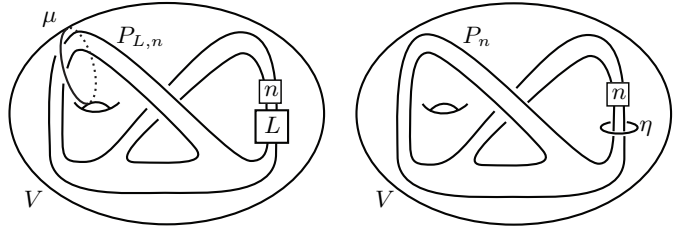

There is a swallow-follow torus tied in the knot $J$; for example, a swallow-follow torus is obtained by tying the solid torus $V$ from Figure~\ref{pattern} into the knot $J$ (and twisting $m$ times).  Figure~\ref{disc} shows the swallow-follow torus that is tied into the knot $L$ (ignore the $D$ for now).  We aim to show that in certain cases, $K$ is a satellite of $J$ (or of $L$).

\begin{prop}\label{sat1}
	Suppose $n=0$.  Then $K$ is a satellite of of $J$ if and only if $L$ is nontrivial.
\end{prop}
\begin{proof}
	If $L$ is the unknot, then $K$  is the unknot, which is not a satellite of $J$.  Now suppose that $L$ is nontrivial.  Consider the swallow-follow torus $V_J$ that contains $K$.  Then $K$ is geometrically essential in $V_J$ if and only if the pattern knot $P$ is geometrically essential in $V$, where $P=P_{L,0}$ as in Figure~\ref{pattern}.  Let $\mu$ denote the meridian of $V$.  Suppose $\mu$ bounds a disc in $V$ missing $P$.  By Dehn's lemma, we may assume this disc $D$ is properly imbedded.  There is a twice punctured disc $N$ bounded by $\mu$ and two meridians of $P$.  Assume that $N$ and $D$ have common boundary and intersect transversely in their interiors. Using a cut-and-paste procedure, we may remove intersection curves that are inessential in $N$.  Let $\beta$ be a curve in $N\cap D$ that is essential in $N$.  Now $\beta$ bounds a subdisc of $D$ that is disjoint from $P$, which means $\mbox{lk}(\beta,P)=0$.

	Claim: This $\beta$ is isotopic in $N$ to $\mu$.  Since $\beta$ is an essential simple closed curve in a planar surface with only three boundary components, $\beta$ must be boundary parallel.  Two of the boundary components algebraically link $P$ in a nontrivial way, so $\beta$ must be parallel to $\mu$.

	Thus, we may (if $\beta$ is outermost in $N$) replace $D$ with a disc that intersects $N$ in fewer essential curves.  We have proven that if $\mu$ bounds a disc in $V-P$, then $\mu$ bounds a disc in $V-(P\cup N)$.  Now cut $V-P$ along $N$ to obtain a ball with two knotted handles burrowed out (see Figure~\ref{split}); call this space $X$.  

\begin{figure}[ht]
  \begin{center}
  	\begin{picture}(350, 100)(0,0)
		\includegraphics{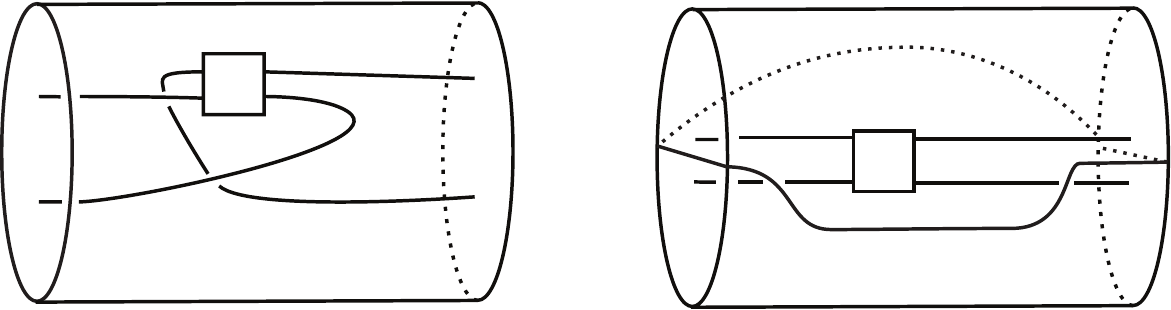}
		\put(-273, 63){$L$}
		\put(-86, 40){$L$}
		\put(-175, 44){\Large $\sim$}
		\put(-346, 41){$\mu$}
		\put(-81, 17){$\mu$}
	\end{picture}
	\caption{The space $X$}
	\label{split}
  \end{center}
\end{figure}
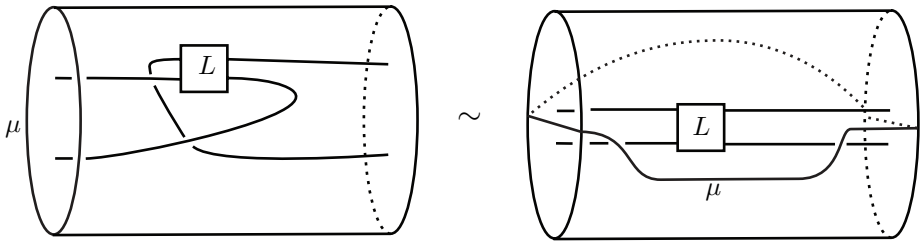

	Since the disc $D$ bounded by $\mu$ missed $N$, $\mu$ bounds a disc in $X$.  Attach two arcs $A$ and $A'$ to the arcs in $X$ to obtain a link that is comprised of two parallel copies of $L$ (see Figure~\ref{splitlink}).  

\begin{figure}[ht]
  \begin{center}
  	\begin{picture}(140, 90)(0,0)
		\includegraphics[scale=.8]{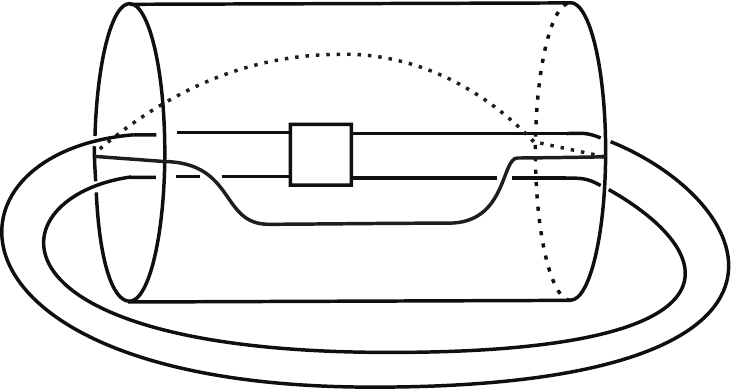}
		\put(-96, 32){\small $\mu$}
		\put(-98, 51){$L$}
	\end{picture}
	\caption{Parallel copies of $L$}
	\label{splitlink}
  \end{center}
\end{figure}

	After viewing Figure~\ref{splitlink}, one sees that the curve $\mu$ bounds a disc $D'$ in $S^3-X-A-A'$, and $D\cup D'$ is a two-sphere $S$.  Since $D$ splits the strands of the link in $X$ and $D'$ splits the strands of the link in $S^3-X$, the two-sphere $S$ splits the link.  Since the link consists of two parallel copies of $L$, there is an annulus $A$ spanning the link.  Since $S$ splits this link, there is a (simple closed) curve $\gamma$ in $A\cap S$ that is essential in $A$, hence isotopic to $L$.  Now $\gamma$ is a knot in a two-sphere, so $\gamma$ and $L$ must be trivial, a contradiction.
\end{proof}

\begin{prop}\label{essential}
	Let $P_n$ and $V$ be as in Figure~\ref{pattern}.  If $n\neq 0$, then $P_n$ is geometrically essential in the solid torus $V$.
\end{prop}
\begin{proof}
	We must show that a meridian $\mu$ of $V$ does not bound a disc in $V-P_n$.  Consider the link $\mu\cup P_n$ pictured below:

\begin{center}
	\begin{picture}(142, 65)(0,0)
		\includegraphics[viewport=0 0 142 57]{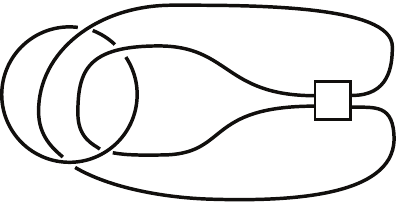}
		\put(-49, 28){\small$n$}
		\put(-146, 41){$\mu$}
		\put(-85, 7){$P_{n}$}
   \end{picture}
\end{center}

	It suffices to show that $\mu$ does not bound a disc in the complement of $P_n$.  By \cite{dR76}, the multivariable Alexander polynomial of this link is $n(1-x)(1-y)$.  For our purposes, it is inconsequential whether $x$ corresponds to the meridian of $\mu$ or of $P_n$.  What is important is that the $x$-degree (and $y$-degree) of this polynomial is one.  Now $\mu$ cannot bound a disc that is embedded in the complement of $P_n$, as Corollary 10.4 of \cite{sH05} implies the first betti number of a surface bounded by $\mu$ will be greater than one.
\end{proof}
  
\begin{lem}\label{alternate}
	Let $P_{L,n}$ and $P_n$ denote the pattern knots in solid torus $V$ pictured in Figure~\ref{pattern}.  The manifold $V-P_{L,n}$ can be obtained by replacing the solid torus $\eta$ in $V-P_n$ with the exterior $E(L)$ of $L$ in such a manner that identifies the longitude of $L$ with the meridian of $\eta$, and the longitude of $\eta$ with the meridian of $L$.
\end{lem}

We ommit the proof of this well-known lemma.
  
\begin{prop}\label{sat2}
	Let $K=K(J,L,m,n)$.  If $n\neq 0$, then $K$ is a satellite of $J$.
\end{prop}
\begin{proof}
	Let us work with the pattern knot $P_{L,n}$ in the solid torus $V$.  We think of the alternate description of $V-P_{L,n}$ from Lemma~\ref{alternate}.  Let $f:V\hookrightarrow S^3$ be an embedding that ties $V$ into the knot $J$ and taking $P_{L,n}$ to $K$ (note that $f$ is faithful if and only if $m=0$).  It suffices to prove $P_{L,n}$ is geometrically essential in $V$, which will be proven by contradiction.

	Assume the meridian $\mu$ of $V$ bounds a disc $D$ in $V-P_{L,n}$ and that $D$ and $\del\eta$ are transverse.  The intersection of $D$ and $\del\eta$ is a collection of simple closed curves.  By a standard cut-and-paste argument, we may eliminate the intersection curves that are inessential in $\del\eta$.  Now we have that $D\cap\del\eta$ consists of simple closed curves $\gamma_1,\ldots,\gamma_l$ that are essential in $\del\eta$.  Since $\del\eta$ is a torus, any two of the $\gamma_i$ cobound exactly two annuli in $\del \eta$.  Since $D$ is a disc, $D-\gamma_i$ is disconnected.  Since $L$ is nontrivial, the component of $D-\gamma_i$ that does not contain $\del D$ must intersect $\del \eta$ somewhere besides $\gamma_i$; thus this component intersects $\del \eta$ in another curve $\gamma_{i'}$.  This proves that there are curves $\gamma_i$ and $\gamma_j$ that cobound an annulus $A\subset D$ with the property that $A$ is properly embedded in $E(L)$.  Let $B\subset \del\eta$ be one of the annuli cobounded by $\gamma_i$ and $\gamma_j$.  Using a product neighborhood of $\del\eta$, let $B'$ denote a small pushoff of $B$ into $V-P_n$; denote by $\gamma_i'$ and $\gamma_j'$ the boundary components of $B'$.  Form an immersed disc $D'$ by gluing $B'$ to the component of $D-\gamma_i'-\gamma_j'$ that contains $\del D$.  Observe that the singularities of $D'$ are contained in $D'-\del D'- (D'\cap\del\eta)$ and that $D'$ intersects $\del\eta$ in $l-2$ curves.

	We may apply this argument until we are left with an immersed disc $D$ with no singularities on the boundary, which is $\mu$, and $D$ is disjoint from $\del\eta$.  Thus, this disc lies in $V-(P_n\cup \eta)$.  By Dehn's lemma, $\mu$ bounds a properly embedded disc in $V-P_n$, contradicting Proposition~\ref{essential}.
\end{proof}

Combinging Propositions~\ref{sat1},~\ref{essential}, and~\ref{sat2}, we have

\begin{thm}\label {satelliteknots}
	Let $K=K(J,L,m,n)$.  Then $K$ is a satellite of $J$ if and only if one of the following hold:
\begin{enumerate}
	\item $n\neq 0$, or
	\item $n=0$ and $L$ is a nontrivial knot.
\end{enumerate}
\end{thm} 

\subsection{Seifert surfaces of $K(J,L,m,n)$}\label{seifsurf}

	There is an obvious genus one Seifert surface for $K$, called the \emph{standard Seifert surface for $K$}.  Theorem~\ref{satelliteknots} gave the necessary and sufficient conditions for $K$ to be a winding number zero satellite of $J$ (and of $L$).  Let $N(J)$ (respectively, $N(L)$) denote the solid companion torus whose core is the knot $J$ (respectively, $L$); unless otherwise specified, we assume that $K$ lies in the interior of $N(J)$.  Observe that $K$ may not lie in both the interior of $N(J)$ and the interior of $N(L)$.  We say $K$ is of \emph{order} 2 with respect to $J$ (respectively, $L$), since there is a properly embedded disc in $N(J)$ (respectively, $N(L)$) that intersects $K$ twice, and 2 is the minimal number of intersection points over all such discs.  Since $K$ is a satellite of $J$, $K$ must be nontrivial and of genus one.  While we do not prove that $K$ has a unique genus one Seifert surface, we aim to say something significant about an arbitrary Seifert surface for $K$.

  In the case $n=\pm1$ and $L$ is the unknot, $K(J,L,m,n)$ is the positively- or negatively-clasped, $m$-twisted Whitehead double of $J$.  In~\cite{wW73}, Whitten showed that Whitehead doubles of $\SR$-knots have a unique genus one Seifert surface up to isotopy.  His first step was to isotope an arbitrary genus one Seifert surface into the neighborhood of the companion.  Here we prove a similar result for $K(J,L,m,n)$ where $J,L\in\SR$.

\begin{thm}\label{movein}
	Let $K=K(J,L,m,n)$ where $J,L\in\SR$.  Then any minimal genus Seifert surface for $K$ may be isotoped to lie in the interior of $V=N(L)$.
\end{thm}
\begin{proof}
	Suppose $K\subset\mathring V$.  Let $S$ be a genus one (minimal genus) Seifert surface for $K$.  Assume $S$ and $T=\del V$ are transverse.  Then $S\cap T$ is a disjoint collection of simple closed curves.  If some intersection curve $\delta$ is inessential in $T$, then $\delta$ bounds a 2-disc $\Delta\subset T$.  Take an innermost (on $T$) such curve $\delta$, so that $\mathring\Delta\cap S=\emptyset$.  By Proposition~\ref{twice}, the fact that $\delta$ is not parallel to $K$ ($\delta$ is the unknot, while $K$ has genus one), $\delta$ bounds a disc $D\subset S$.  The 2-sphere $\Delta\cup D$ bounds a 3-cell disjoint from $S-D$, since $\delta$ is innermost.  We may use this 3-cell to isotope $D$ across $T$, hence eliminating the curve $\delta$ from the intersection of $S$ and $T$.

	Having eliminated all inessential intersection curves by isotoping $S$, we see $S\cap T$ is a disjoint collection of curves $\gamma_1,\ldots,\gamma_k$, each of which is essential on the torus $T$.  We may order them so that for each $i$, $\gamma_i\coprod\gamma_{i+1}$ bounds an annulus in $T$ that does not intersect the other $\gamma$ curves.  We see that these $\gamma_i$ are parallel on the torus $T$, and we will see presently that they are also parallel on the punctured torus $S$.  The two claims below imply that each $\gamma_i$ is essential in $H_1(S)$.  By the classification of the curve complex of the punctured torus~\cite{yM99}, the $\gamma_i$ must be isotopic (parallel) in $S$.

	Claim 1: No $\gamma_i$ bounds a disc in $S$.  If $\gamma_i$ did bound a disc in $S$, the disc could not lie inside $V$, as any properly embedded disc in $V$ would intersect $K$, since $K$ is a satellite of $L$.  Therefore this disc would be properly embedded in $S^3-\mathring V$, implying $\gamma_i$ is a longitude of $L$ and that $L$ has genus zero, contradicting the assumption $L$ is nontrivial.

	Claim 2: No $\gamma_i$ is parallel in $S$ to $K$.  This follows from the fact that $\gamma_i$ is an embedded, essential curve in $T$ and $K$ has winding number zero in $V$ and is geometrically essential in $V$.  Thus, $K$ cannot be isotoped in $V$ to $\gamma_i$.

	Claims 1 and 2 imply each $\gamma_i$ is essential in $H_1(S)$, so $S-\gamma_i$ is connected.  Since the $\gamma_i$ are parallel on $S$ and essential in $H_1(S)$, $S-\cup \gamma_i$ is a twice punctured disc together with a collection of disjoint annuli (there are, in fact, $k-1$ of these annuli).  As per Whitten's proof, there are three cases:
	\begin{itemize}
		\item[(i)] each $\gamma_i$ has winding number 1 in $V$,
		\item[(ii)] each $\gamma_i$ is a cable of $L$, or
		\item[(iii)] each $\gamma_i$ is a meridian of $V$.
	\end{itemize}

	Whitten's proof for cases (i) and (ii) works in our current situation, so we omit these cases and refer the reader to \cite{wW73}.  We show presently that case (iii) cannot occur.

	If $\gamma_1$ is a meridian of $V$, let $D$ be the disc in $V$ bounded by $\gamma_1$ that intersects $K$ in 2 points (see Figure~\ref{disc}).  Let $\alpha$ be a small pushoff of $\gamma_1$ that lies outside of $V$ and does not intersect $S$.  Then $\alpha$ bounds a disc $D'$ with $D\subset D'$.  We may assume $S$ and $D'$ are transverse, so  $S\cap D'$ is a collection of simple closed curves and one arc $\sigma$, the boundary of which is two points on $K$:
	\begin{center}
		\begin{picture}(130, 70)(7, 2)
			\put(41, 33){$\sigma$}
			\put(46, 48){\small $\zeta$}
			\includegraphics[scale=.8]{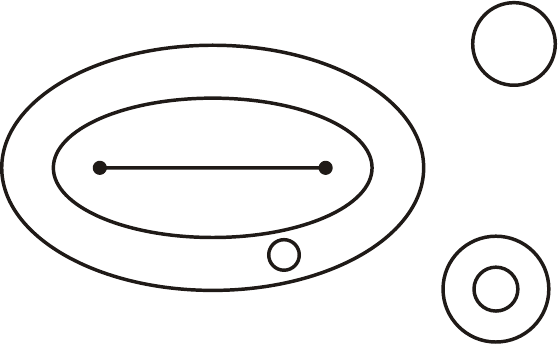}
		\end{picture}
	\end{center}

	Suppose $\delta$ is one of these simple closed curves in $S\cap D'$.  Assume $\delta$ is an innermost curve with the property that the subdisc $\Delta$ of $D'$ bounded by $\delta$ misses $\sigma$.  This $\delta$ must also bound a disc $\Delta'\subset S$, and $\Delta\cup\Delta'$ bounds a 3-cell that intersects $S$ only in $\Delta'$.  We may isotope $S$ through this 3-cell to eliminate the curve $\delta$ from the intersection of $S$ and $D'$.

	\begin{figure}[ht]
		\begin{center}
			\begin{picture}(140, 135)(0,0)
				\includegraphics[viewport=0 0 234 216, scale=.5]{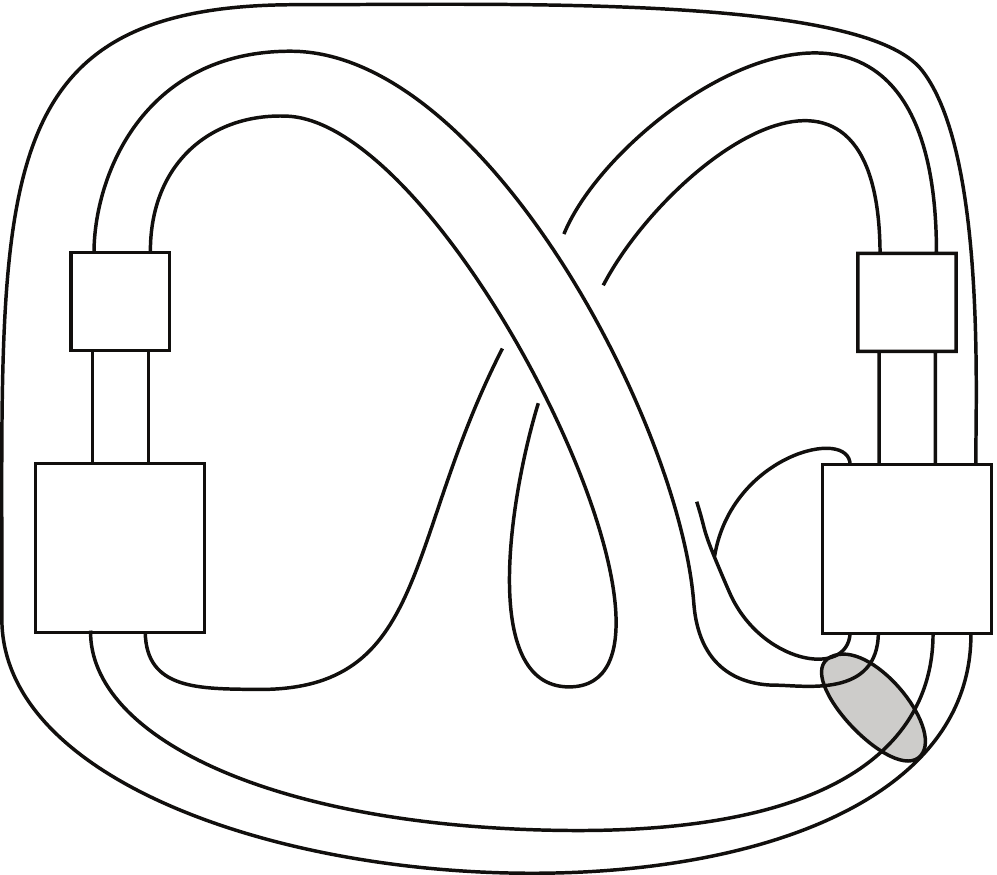}
				\put(-102, 44){$J$}
				\put(-104, 81){$m$}
				\put(10, 45){$L$}
				\put(11, 81){$n$}
				\put(15, 10){$D$}
				\put(-134, 14){$N(L)$}
			\end{picture}
			\caption{The knot $K$ in the companion torus $N(L)$ with disc $D$}
			\label{disc}
		\end{center}
	\end{figure}
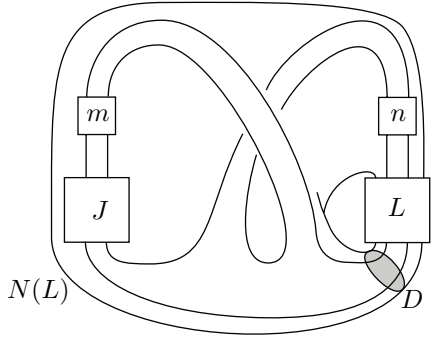

	Now we may assume that $S\cap D'$ is the arc $\sigma$ together with some closed curves that encircle $\sigma$.  Choose an innermost such curve $\zeta$.  Now $\zeta$ is not parallel in $S$ to $K$, since $\zeta$ is the unknot. Since the unknot is fibered, any disc bounding $\zeta$ may be isotoped inside $V$.  Since $K$ is geometrically essential in $V$, any disc in $V$ bounding $\zeta$ must intersect $K$, and so $\zeta$ does not bound a disc in $S$.  We have shown that $\zeta$ does not separate $S$.  Cut $S$ along the disc $\Delta\subset D'$ with $\del \Delta=\zeta$; this cuts $S$ along $\sigma$ (see Figure~\ref{cutandglue}).  Glue two copies of $\Delta$ and two copies of $\sigma$ to $S-\Delta$ to obtain a (disconnected) surface $\Sigma$ spanning the link that consists of two parallel copies of $J$ (see Figure~\ref{parallels}); here our choice of $D$ is crucial.

	\begin{figure}[ht]
		\begin{center}
			\includegraphics[scale=.5]{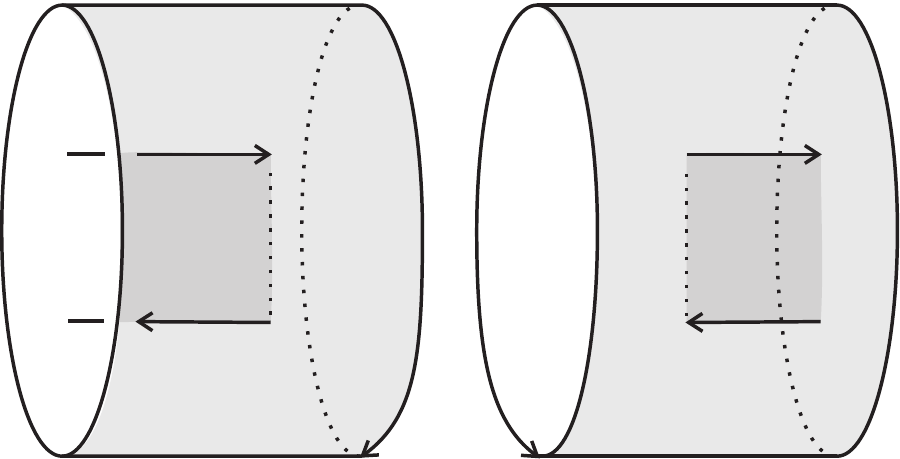}
		\end{center}
		\caption{$S$ cut along the disc $\Delta$}
		\label{cutandglue}
	\end{figure}

	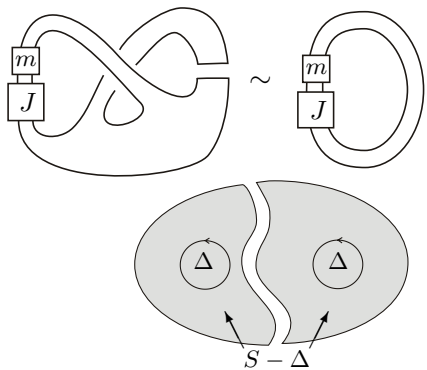
\begin{figure}[ht]
		\begin{center}
			\begin{picture}(220, 100)(0,0)
				\includegraphics[viewport=0 0 288 144, scale=.6]{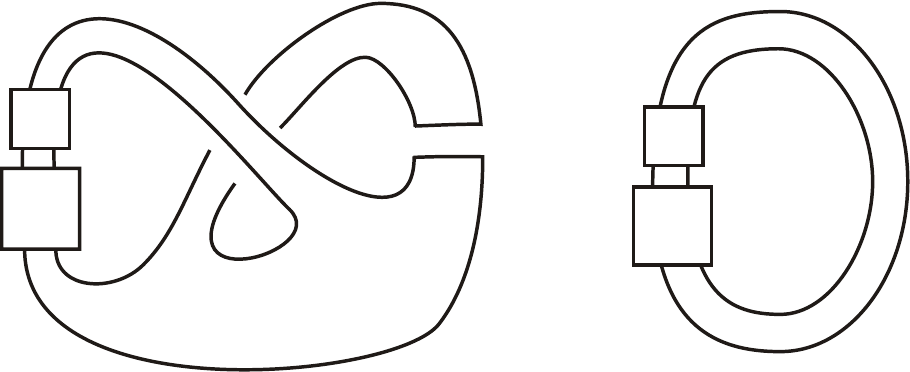}
				\put(-169, 25){$J$}
				\put(-170.5, 42){\small $m$}
				\put(-59, 22){$J$}
				\put(-60.5, 39){\small $m$}
				\put(-82, 34){\large $\sim$}
			\end{picture}
		\end{center}
		\begin{center}
			\begin{picture}(110, 50)(5,5)
				\includegraphics[viewport=5 5 265 191, scale=.4]{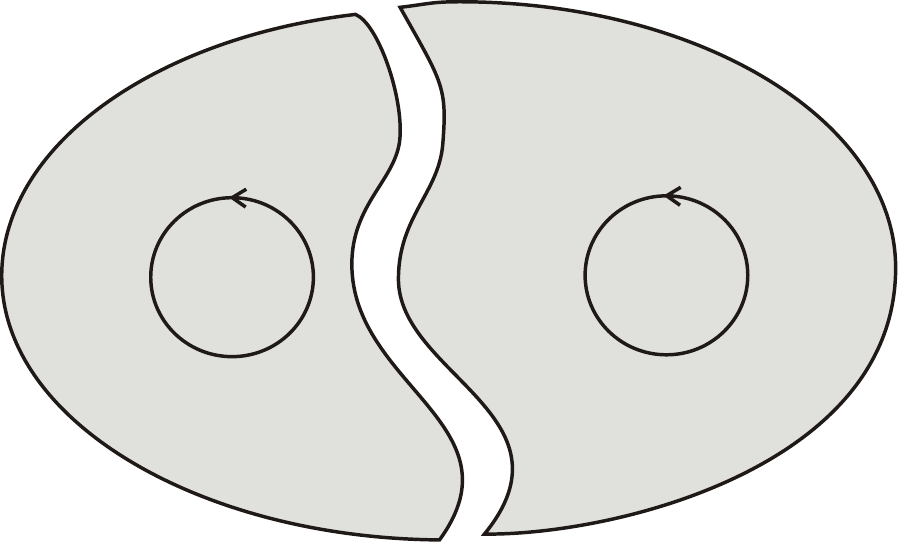}
				\put(-84, 27){\small $\Delta$}
				\put(-33, 27){\small $\Delta$}
				\put(-65, -10){\small $S-\Delta$}
				\put(-65, -3){\vector(-1,2){7}}
				\put(-40, -3){\vector(1,2){7}}
			\end{picture}
		\end{center}
		\caption{The boundary of the new surface $\Sigma$ and a schematic of $\Sigma$}
		\label{parallels}
	\end{figure}

	Since $S$ has genus one, $\Sigma$ must be the disjoint union of two discs, hence $J$ is the unknot, contradicting the nontriviality of $J$.  We conclude that $\zeta$ cannot be a meridian of $V$.
\end{proof}

\begin{defn}
	Let $A$ be an oriented annulus whose core curve is a knot $J$ and such that the boundary curves have linking number $m$.  Define the oriented two-component link $M_J^m$ to be the boundary of $A$.  Figure~\ref{parallels} depicts $M_J^m$.  There is an obvious annulus in Figure~\ref{parallels} spanning $M_J^m$; call this annulus the \textbf{standard annulus} spanning $M_J^m$.
\end{defn}

\begin{lem}\label{arcs}
	Let $A$ and $B$ be arcs in $\R^2$ with common boundary, and assume $A$ and $B$ intersect transversely.  Enumerate the intersection points $a_0,\ldots,a_n\in A$ and $b_0,\ldots,b_n\in  B$ so that
	\begin{itemize}
		\item $a_0=b_0$ and $a_n=b_n$,
		\item $\del A=\del B=a_0\cup a_n=b_0\cup b_n$, and
		\item for each $i$, $a_i$ and $a_{i+1}$ are adjacent as points on $A$, and $b_i$ and $b_{i+1}$ are adjacent as points on $B$.
	\end{itemize}

	If $\sigma$ denotes the permutation determined by $a_{i}=b_{\sigma(i)}$ for each $i=0,\ldots,n$, then for some $i$, $\left|\sigma(i+1)-\sigma(i)\right|=1$; that is, $b_{\sigma(i)}$ and $b_{\sigma(i+1)}$ are adjacent on $B$.
\end{lem}
\begin{proof}
	Let $\sigma$ denote the permutation described above.  First note that $\sigma(0)=0$ and $\sigma(n)=n$.  The claim is clearly true if $n=1$, so assume $n>1$.  Let $A_{1}$ denote the subarc of $A$ bounded by $a_{0}$ and $a_{\sigma\inv(1)}$, and let $B_{1}$ denote the subarc of $B$ bounded by $b_{0}$ and $b_{1}$.  Now $A_{1}\cup B_{1}$ is a circle that bounds precisely one disc $D_{1}$ in $\R^{2}$.  Either $a_{n}\in D_{1}$ or $a_{n}\not\in D_{1}$.
	
	Assume $a_{n}\in D_{1}$.  The reader may refer to Figure~\ref{arcdiagram} for a picture.  Let $A_{n}$ denote the subarc of $A$ bounded by $a_{\sigma\inv(n-1)}$ and $a_{n}$, and let $B_{n}$ denote the subarc of $B$ bounded by $b_{n-1}$ and $b_{n}$.  Since $B_{n}$ contains none of the $b_{i}$ in its interior, $B_{n}$ must lie in $D_{1}$.  Since $a_{n}\in D_{1}$, the entire arc $A$ must lie in $D_{1}$.  In particular, $A_{n}\subset D_{1}$.  Let $D_{n}$ denote the disc in $\R^{2}$ bounded by the circle $A_{n}\cup B_{n}$.  Observe that $D_{n}\subset D_{1}$ and $a_{0}\not\in D_{n}$.
	
	\begin{figure}[ht]
		\begin{center}
			\begin{picture}(255, 85)(0,0)
				\includegraphics{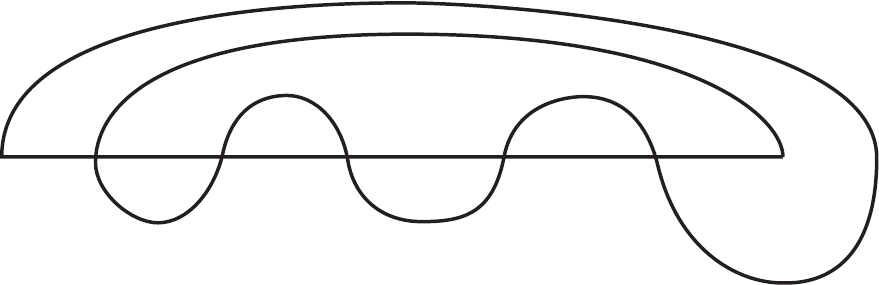}
				\put(-260, 31){$a_{0}$}
				\put(-32, 31){$a_{n}$}
				\put(-110, 8){$b_{1}=a_{\sigma\inv(1)}$}
				\put(-282, 8){$b_{n-1}=a_{\sigma\inv(n-1)}$}
				\put(-235, 16){\vector(1,3){6}}
				\put(-72, 16){\vector(1,3){6}}
				\put(-32, 8){$B_{1}$}
				\put(-139, 64){$B_{n}$}
			\end{picture}
			\caption{The arcs $A$ and $B$.}
			\label{arcdiagram}
		\end{center}
	\end{figure}
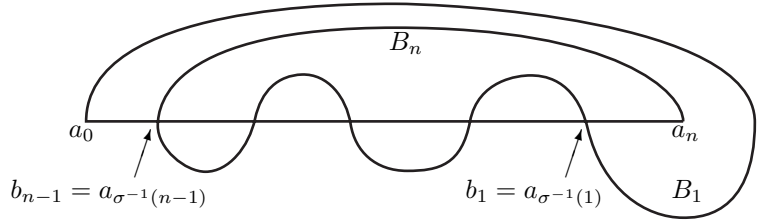
	
	We have established that $a_{0}\not\in D_{n}$ or $a_{n}\not\in D_{1}$.  Without loss of generality, assume $a_{0}\not\in D_{n}$.
	
	Now $\left|\sigma\inv(n-1)-\sigma\inv(n)\right|=\left|\sigma\inv(n-1)-n\right|>0$.  If $\left|\sigma\inv(n-1)-\sigma\inv(n)\right|=1$, we may proceed with the next paragraph of the proof.  If $\left|\sigma\inv(n-1)-\sigma\inv(n)\right|>1$, there must be some $a_{i}$ that lie in the interior of $A_{n}$.  Thus, there is a subarc $B'$ of $B$ that is properly embedded in $D_{n}$ with boundary in the interior of $A_{n}$.  

Write $\del B'=a_{\sigma\inv(j)}\cup a_{\sigma\inv(j-1)}$ for some $j$.  Observe that $$\left|\sigma\inv(j-1)-\sigma\inv(j)\right|<\left|\sigma\inv(n-1)-\sigma\inv(n)\right|$$
	
	By the preceding paragraph, we can find a $k$ such that $\left|\sigma\inv(k-1)-\sigma\inv(k)\right|=1$, and the proof is complete.
\end{proof}

\begin{lem}\label{annulus}
	If $J\in\SR$, then any two annuli with common boundary $M_J^m$ are isotopic rel boundary.
\end{lem}
\begin{proof}
	Let $A$ and $B$ be annuli with common boundary $M_J^m$.  Assume the interiors of $A$ and $B$ intersect transversely, so that $A\cap B$ is a collection of disjoint simple closed curves on $A$, including the two boundary components.  If any intersection curve bounds a disc in $A$, it must also bound a disc in $B$.  Picking an innermost curve on $A$, we may isotope $B$ to eliminate the curve from the intersection without creating new intersection curves.
    
	We now have $A\cap B=a_0\coprod \cdots\coprod a_n$, where $a_0\cup a_n=\del A=\del B$.  Assume that we have ordered the $a_i$ so that $a_i\cup a_{i+1}$ is the boundary of a subannulus $A_{i+1}$ of $A$ with the property that $A_{i+1}\cap a_j=\emptyset$ if $j\neq i,i+1$.  We say that $a_i$ and $a_{i+1}$ are \emph{adjacent} on $A$, and we call $A_{i+1}$ an \emph{adjacency annulus} on $A$.
    
	If $n=1$, then $A\cup B$ is an embedded torus, which must bound a solid torus $V$.  Since $J$ is isotopic to $a_0$ and is not a cable, $J$ must be a winding number one satellite of the core of $V$.  Write $[J]=[\lambda]+q[\mu]\in H_1(\del V)$, where $\lambda$ is a core curve of $V$, and $\mu$ is a meridian of $\lambda$.  By a theorem of Schubert~\cite[Satz, p. 165]{hS53}, $A$ is isotopic rel boundary to $B$ in $V$, completing the proof for the case $n=1$.
    
	Assume $n>1$.  Let $c$ be a nonseparating arc that is properly embedded in $A$.  Pick a nonseparating arc $d$ that is properly embedded in $B$ so that $d$ has the same boundary as $c$ and so that $c\cup d$ may be isotoped to lie on a disc.  We use this disc to take a cross-section of $A\cup B$.  This cross-section is the union of two arcs which have common boundary and which are transverse.  By Lemma~\ref{arcs}, we may find a pair of points in the intersection of $c$ and $d$ that are adjacent on each of $c$ and $d$.  Thus, in the context of the present argument, we may find a pair of curves in the intersection of $A$ and $B$ that are adjacent on each of $A$ and $B$.  Denote these curves as $a_{i-1}$ and $a_{i}$, and denote the adjacency annuli $A_{i}$ and $B_{i}$.  Let $V_i$ be solid torus bounded by $A_i\cup B_i$.  Observe that $V_{i}\cap (A\cup B)=A_{i}\cup B_{i}$.   We argued previously that $J$ is a winding number 1 satellite of $V_i$, thus $V_i$ is the unique solid torus with boundary $A_i\cup B_i$.  Furthermore, the core of $V_i$ is isotopic to $J$.
    
	By Schubert's theorem, we may isotope $A_i$ to $B_i$ while fixing $a_{i-1}$ and $a_i$ pointwise.  If $i\neq n$ and $i\neq 1$, we must push $A_i$ through $B_i$ to maintain the transversality of $A$ and $B$.  If $i=n$ (or $i=1$), we must push $A_n-a_n$ through $B_n-a_n$ while leaving $a_0$ fixed.  Either case may be easily accomplished; Figure~\ref{easy} depicts cross-sections of how one might go about accomplishing the task.  During this isotopy, $A_i$ never intersects any of the other $A_j$, since $A_i$ is an adjacency annulus.  Thus, this isotopy of $A_i$ extends to an isotopy of $A$, fixing the boundary.  The isotopy reduces the number of intersection curves, and the proof is complete.

	\begin{figure}[ht]
		\begin{center}
			\begin{picture}(252, 100)(0,0)
				\put(35, 84){\large $A_i$}
				\put(23,62){\large $B_i$}
				\put(15,31){\large $A_n$}
				\put(29, 8){\large$B_n$}
				\put(52, 10){\large $a_n$}
				\includegraphics{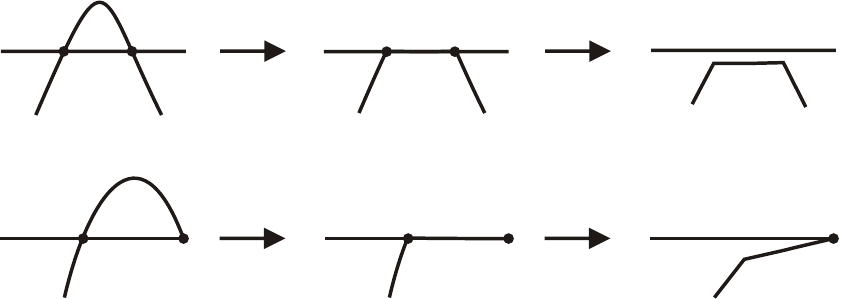}
			\end{picture}
		\end{center}
		\caption{Pushing $A_i$ through $B_i$}
		\label{easy}
	\end{figure}
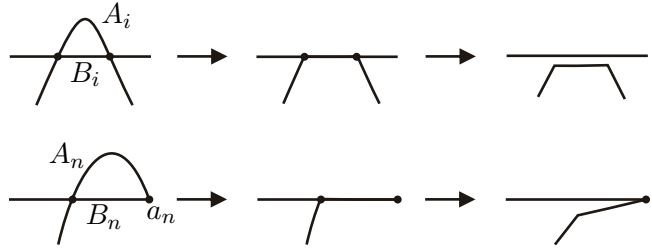
\end{proof}

\begin{thm}\label{sbtheorem}
	If $J,L\in\SR$ and $S$ is a genus one Seifert surface for $K=K(J,L,m,n)$, then $S$ has a symplectic basis $(\alpha,\beta)$ where $\alpha$ and $J$ have the same knot type, and $\beta$ and $L$ have the same knot type.
\end{thm}
\begin{proof}
	By Theorem~\ref{movein}, we may assume $S\subset \mathring N(L)$. As in the proof of Theorem~\ref{movein}, there is a disc $D$ properly embedded in $N(L)$ with $D\cap S=\sigma$, an arc properly embedded in $S$, and $S-\left(\sigma\times(-\epsilon,\epsilon)\right)$ is an annulus $B$ bounding $M_J^m$ (see Figure~\ref{parallels}).  Lemma~\ref{annulus} permits us to assume $B$ is the standard annulus spanning $M_J^m$.  Let $\beta$ be the core curve of $B$, so that $\beta\sim J$.  Pick a point $x\in \mathring \sigma$.  Again using Lemma~\ref{annulus}, let $\alpha'$ denote a  properly embedded arc in $B$ joining $x\times-\epsilon$ to $x\times \epsilon$ that intersects $\beta$ in only one point.  Let $\alpha=\alpha'\cup\left(x\times(-\epsilon,\epsilon)\right)$.  Since $S$ may be isotoped into $N(J)$, we may pick $\alpha'$ so that $\alpha$ is isotopically trivial (i.e. geometrically inessential) in $N(J)$.  From a geometric point of view, we must carefully pick $\alpha'$ so that $\alpha$ ``does not wrap around $\beta$."  Observe $\alpha$ is a simple closed curve intersecting $\beta$ in one point.  We conclude $(\alpha,\beta)$ is a symplectic basis for $S$.

	Now isotope $S$ to lie in $\mathring N(J)$.  Since $\alpha$ is isotopically trivial in $N(J)$, there is a disc $D$ properly embedded in $N(J)$ such that $D\cap S=\delta$ is an arc properly embedded in $S$ with $\delta\cap\alpha=\emptyset$ and $\delta\cap\beta$ is a point.  Cutting $S$ along $D$ yields an annulus $A$ bounding $M_L^n$, and the curve $\alpha$ lies on $A$.  By Lemma~\ref{annulus}, $\alpha\sim L$.
\end{proof}

\begin{cor}\label{bound}
	For $J,L\in\SR$, $g_1(K(J,L,m,n))\geq g(J)+g(L)$.
\end{cor}
\begin{proof}
	Let $S$ be a genus one Seifert surface for $K=K(J,L,m,n)$, and equip $S$ with a symplectic basis $(\alpha,\beta)$ from Theorem~\ref{sbtheorem}.  We may write any basis for $H_1(S)$ as $x=p\alpha+q\beta$ and $y=r\alpha+s\beta$, where $$\det\left(\begin{array}{rr}
      p&q\\
      r&s
\end{array}\right)=1$$

	Without loss of generality, assume $p\neq 0$.  Since $S$ may be isotoped to lie in $N(J)$, $x$ is a winding number $p\neq 0$ satellite of $J$.  Similarly, if $s\neq 0$, then $y$ is a winding number $s\neq 0$ satellite of $L$.  On the other hand, if $s=0$, then $r\neq 0$ and $q\neq 0$, in which case $y$ is a winding number $r\neq 0$ satellite of $J$, and $x$ is a winding number $q\neq 0$ satellite of $L$.  In either case, $g(x)+g(y)\geq g(J)+g(L)$.
\end{proof}

\begin{exam}
	Let $K=K(J,L)$ where $J,L\in\SR$.  We claim $g_1(K)=g(J)+g(L)$.  Let $S$ be the standard Seifert surface for $K$.  We may construct a weak grope for $K$ by attaching a minimal genus Seifert surface for $J$ to $S$ from the front, and a surface for $L$ to $S$ from the back.  The sum of the genera of these two surfaces is $g(J)+g(L)$, and  Corollary~\ref{bound} allows us to conclude that $g_1(K)=g(J)+g(L)$.
\end{exam}

%%%%%%%%%%%%%%%%%%%%   End of main body of article
%
%                             References
%
%   BiBTeX users uncomment the following line:
%
%\bibliographystyle{gtart}
%

%\begin{thebibliography}
	\bibliographystyle{amsalpha}
	\bibliography{fog}

\providecommand{\bysame}{\leavevmode\hbox to3em{\hrulefill}\thinspace}
\providecommand{\MR}{\relax\ifhmode\unskip\space\fi MR }
\providecommand{\MRhref}[2]{%
  \href{http://www.ams.org/mathscinet-getitem?mr=#1}{#2}
}
\providecommand{\href}[2]{#2}
\begin{thebibliography}{COT03}

\bibitem[BJ06]{BJ06}
Mark Brittenham and Jacqueline Jensen, \emph{Canonical genus and whitehead
  doubles of pretzel knots}, Available at http://arxiv.org/abs/math/0608765,
  2006.

\bibitem[Coc04]{tC04}
Tim Cochran, \emph{Noncommutative knot theory}, Algebraic \& Geometric Topology
  \textbf{4} (2004), 347--398.

\bibitem[COT03]{COT03}
Tim Cochran, Kent Orr, and Peter Teichner, \emph{Knot concordance, whitney
  towers and $l^2$-signatures}, Annals of Mathematics \textbf{157} (2003),
  433--519.

\bibitem[CT07]{CT07}
Tim Cochran and Peter Teichner, \emph{Knot concordance and von neumann
  $rho$-invariants}, Duke Mathematical Journal \textbf{137} (2007), no.~2,
  337--379.

\bibitem[Har05]{sH05}
Shelly Harvey, \emph{Higher-order polynomial invariants of 3-manifolds giving
  lower bounds for the thurston norm}, Topology \textbf{44} (2005), 895--945.

\bibitem[Hem76]{jH76}
John Hempel, \emph{3-manifolds}, Annals of Mathematics Studies, no.~86,
  Princeton University Press, 1976.

\bibitem[Kaw84]{aK84}
Akio Kawauchi, \emph{Classification of pretzel knots}, Kobe Journal of
  Mathematics \textbf{2} (1984), 11--22.

\bibitem[Min99]{yM99}
Yair Minsky, \emph{The classification of punctured-torus groups}, Annals of
  Mathematics \textbf{149} (1999), 559--626.

\bibitem[Rol76]{dR76}
Dale Rolfsen, \emph{Knots and links}, Publish or Perish, Berkeley, CA, 1976.

\bibitem[Sch53]{hS53}
Horst Schubert, \emph{Knoten und vollringe}, Acta Mathematica \textbf{90}
  (1953), 131--286.

\bibitem[Tei04]{pT04}
Peter Teichner, \emph{What is a grope?}, Notices of the American Mathematical
  Society \textbf{51} (2004), no.~8, 892--893.

\bibitem[Whi73]{wW73}
Wilbur Whitten, \emph{Isotopy types of knot spanning surfaces}, Topology
  \textbf{12} (1973), 373--380.

\end{thebibliography}
%\end{thebibliography}

\end{document}